\documentclass[12pt]{amsart}

\usepackage{hyperref}
\usepackage{amsmath, amssymb, dsfont, mathabx}

\usepackage{amsthm}
\usepackage{amstext}
\usepackage{amssymb}
\usepackage{mathrsfs}
\usepackage{mathabx}
\usepackage[english]{babel}
\usepackage{tikz-cd}
\usepackage{enumerate}
\usepackage{xfrac}
\usepackage{mathabx}
\usepackage{mathtools}
\usepackage{xcolor}
\usepackage[all]{xy}

\newtheorem*{theorem*}{Theorem}
\newtheorem{theorem}{Theorem}

\newtheorem{definition}{Definition}

\newcommand{\R}{\ensuremath{\mathbb{R}}}
\newcommand{\E}{\ensuremath{\mathbb{E}}}
\DeclareMathOperator{\partr}{/\!/}
\DeclareMathOperator{\ric}{Ric}
\DeclareMathOperator{\End}{End}

\DeclareMathOperator{\id}{id}

\DeclareMathOperator*{\stlim}{st-lim}

\begin{document}

\setlength{\parindent}{0pt}
\title{A note on the scattering theory of Kato-Ricci manifolds}
\author{Batu Güneysu}
\author{Maxime Marot}

\address{Fakultät für Mathematik, Technische Universität Chemnitz, 09126 Chemnitz, Germany}
\email{batu.guneysu@mathematik.tu-chemnitz.de}
\email{maxime.marot@mathematik.tu-chemnitz.de}

\begin{abstract}
In this note we prove a new $L^1$ criterion for the existence and completeness of the wave operators corresponding to the Laplace-Beltrami operators corresponding to two Riemannian metrics on a fixed noncompact manifold. Our result relies on recent estimates on the heat semigroup and its derivative, that are valid if the negative part of the Ricci curvature is in the Kato class - so called Kato-Ricci manifolds.
\end{abstract}

\maketitle

Let $M$ be a connected noncompact manifold of dimension $m\geq 2$. Given a Riemannian metric $g$ on $M$ let $\mu_g$ denote its volume measure, given locally in terms of $\mathbf{g}_{jk}:=g(\partial_j,\partial_k)$ by $d\mu_g(x)=  \sqrt{\det(\mathbf{g})} dx$. The shorthand $\mu_g(x,r)$ denotes the measure of the open ball with respect to the geodesic distance $d_g$ centered in $x$ of radius $r$. The Laplace-Beltrami operator $-\Delta_g$ is a symmetric nonnegative operator in the complex Hilbert space $L^2(M,g):=L^2(M,\mu_g)$, defined initially on $C^\infty_c(M)$. Thus $-\Delta_g$ has a canonically given self-adjoint extension - the Friedrichs realization - which will be denoted with the same symbol again. Then $P_t:=e^{t\Delta_g}$, the right hand side being given via spectral calculus, is the minimal heat semigroup in $L^2(M,g)$, and $P_g(t,x,y)\geq 0$ will denote its integral kernel, so that 
$$
P^g_tf(x)=\int_M f(y)P_g(t,x,y)  d\mu_g(y)\quad  \text{for all $t>0$, $f\in L^2(M,g)$, $x\in M$}.
$$
With the gradient operator $\nabla_g$, we will also be interested in the bounded operator from $L^2(M,g)$ to the Hilbert space of $L^2$-vector fields on $(M,g)$, given by $\hat{P}_t^g:=\nabla_g P_t^g$. Its smooth integral kernel will be denoted with $\hat{P}_t^g(x,y)\in T_x M$, so that
$$
\hat{P}^g_tf(x)=\int_M f(y)\hat{P}_g(t,x,y) d\mu_g(y)\quad \text{for all $t>0$, $f\in L^2(M,g)$, $x\in M$}.
$$
\begin{definition} With 
\begin{align*}
&\sigma_g:M\longrightarrow \R,\quad \sigma_g(x) := \text{\emph{smallest eigenvalue of} } \ric_g(x),\\
&\sigma^-_g:M\longrightarrow [0,\infty),\quad \sigma^-_g(x):=-\min(\sigma_g(x),0),
\end{align*}
we say that \emph{$(M,g)$ is a Kato-Ricci manifold}, if 
\begin{align}\label{norm}
\lim_{t\to 0+} \sup_{x\in M} \int_0^t \int_M P_g(t,x,y) \sigma^-_g(y)\,d\mu_g(y)dt = 0.
\end{align}
\end{definition}

Note that the Ricci curvature of $(M,g)$ is bounded from below by a constant, if and only if $\sigma^-_g$ is bounded, which in view of 
\begin{align}\label{assq}
\int_M P_g(t,x,y)d\mu(y)\leq 1\quad\text{for all $x\in M$, $t>0$,}
\end{align}
implies the Kato-Ricci property. In this sense, Kato-Ricci manifolds generalize manifolds with Ricci curvature bounded from below by a constant. On the other hand, Kato-Ricci manifolds appear quite naturally: for example, one can take a conformal perturbation of the Euclidean metric such that the Laplacian of the conformal factor has higher and higher oscillations on smaller and smaller regions (cf. Example XIV.27 in \cite{bookbatu}).\vspace{1mm}

This note is a contribution at the level of scattering theory to the recent agenda which aims at showing that many results that hold for manifolds with Ricci curvature bounded from below by a constant also hold for Kato-Ricci manifolds. Most importantly, this is confirmed that by the recent breakthrough from \cite{carron_kato_2023}, which states that every complete Kato-Ricci manifold is quasi-isometric in the sense of weighted Riemannian manifolds to a weighted complete Riemannian manifold which, for some (large) $N>1$, has a weighted $N$-dimensional Ricci curvature bounded from below. Other important results of this agenda are:
\begin{itemize}
\item on Kato-Ricci manifolds one has the heat kernel characterzation of the Riemannian total variation \cite{pallara},
\item Kato-Ricci manifolds are stochastically complete \cite{braun_heat_2021}, that is, one has equality in (\ref{assq}),
\item renormalized pointed and measured Gromov-Hausdorff limits of $m$-dimensional Kato-Ricci manifolds with a uniform control on (\ref{norm}) are $k$-rectifiable for some $k\leq m$, and have almost everywhere uniquely determined metric measure tangent cones given by the $k$-dimensional Euclidean space \cite{carron_kato_2023}.    
\end{itemize}
In addition, applications of Kato-Ricci bounds to topology are suggested in \cite{stollmann} and the study of abstract spaces that satisfy the Kato-Ricci assumptions has been initiated in \cite{sturm}.\vspace{1mm}

In order to formulate our main result, assume we are given two complete Riemannian metrics $g,h$ on $M$ and a bounded operator 
$$
J:L^2(M,g)\longrightarrow  L^2(M,h).
$$
Then the wave operators with respect to the identification operator $J$ are defined as the strong limits
$$
W_\pm(-\Delta_h,-\Delta_g;J):=\stlim_{t\to \pm\infty} \exp(-it\Delta_h)J\exp(it\Delta_g)\pi_g^{\mathrm{ac}},
$$
in case these limits exists. Note that, unlike classical scattering problems (that is, scattering by potentials), this scattering problem takes place in two Hilbert spaces. Here for $l=g,h$ the operator $\pi_l^{\mathrm{ac}}$ denotes the projection onto the closed subspace $L^2_{\mathrm{ac}}(M,l)\subset L^2(M,l)$ given by all $f\in L^2(M,l)$ such that, with $E_l$ the projection valued spectral measure of $-\Delta_l$, the Borel measure 
$$
\int_M\left|E_l(A)f(x)\right|^2 d\mu_l(x),\quad A\subset \mathbb{R},
$$
is absolutely continuous with respect to the Lebesgue measure. Given their existence, one says that the wave operators $W_\pm(-\Delta_h,-\Delta_g;J)$ are \emph{complete}, if
$$
\mathrm{Ker}\big(W_\pm(-\Delta_h,-\Delta_g;J))^\perp=L^2_{\mathrm{ac}}(M,g),\quad \overline{\mathrm{Ran}\big(W_\pm(-\Delta_h,-\Delta_g;J))}=L^2_{\mathrm{ac}}(M,h).
$$
The completeness of the wave operators implies \cite{kato} the equality of the two underlying absolutely continuous spectra,
\begin{align}\label{qwcy}
\mathrm{spec}_{\mathrm{ac}}(-\Delta_h)=\mathrm{spec}_{\mathrm{ac}}(-\Delta_g),
\end{align}
where by definition the absolutely continuous spectrum of $-\Delta_l$ is the spectrum of its restriction to $L^2_{\mathrm{ac}}(M,l)$. The physical importance of the absolutely continuous spectrum stems from the fact that absolutely continuous states are scattering states in the physical sense, meaning that if $\psi\in L^2_{\mathrm{ac}}(M,l)$, then 
$$
\lim_{t\to\infty}\int_K|\exp(it\Delta_l)\psi(x)|^2 d\mu(x)=0\quad\text{for all compact $K\subset M$}.
$$
The latter means that the underlying quantum particle eventually leaves every compact. As for a general Riemannian metric standard spectral techniques such as the Fourier transform are not available, results such as (\ref{qwcy}) are very important, as they often allow to determine the absolutely continuous spectrum nevertheless: namely, one can often compare the given metric to one that actually admits Fourier analysis techniques. \vspace{1mm}

In order to formulate our criterion for the existence/completeness of the wave operators, given the locally defined symmetric matrix-valued functions
$$
\mathbf{g}_{jk}=g(\partial_j,\partial_k),\quad \mathbf{h}_{jk}=h(\partial_j,\partial_k),
$$
let $\mathscr{A}_{g,h}$ be the globally defined smooth $(1,1)$ tensor field on $M$ which is locally defined by the positive definite matrix-valued function $\mathbf{g}\mathbf{h}^{-1}$. The fiberwise given endomorphism
$$
\mathscr{A}_{g,h}(x):T_xM\longrightarrow  T_x M
$$
induces the Borel measurable function
$$
\delta_{g,h}:M\longrightarrow [0,\infty), \quad \delta_{g,h}(x):=2\sinh\Big((m/4)\max_{\lambda\in \mathrm{spec}(\mathscr{A}_{g,h}(x))}|\log(\lambda)|\Big),
$$
which measures the deviation of the metrics from each other in a multiplicative and zeroth order way. We write $d\mu_h= \rho_{g,h} \,d\mu_g$, where the smooth density function $\rho_{g,h}:M\to (0,\infty)$ is defined via 
$$
\rho_{g,h}:=\sqrt{\det(\mathbf{h} )   }/\sqrt{\det(\mathbf{g})}.
$$
If $g$ is quasi-isometric to $h$, that is, there exists a constant $C\geq 1$ such that 
$$
(1/C) h \leq g \leq C h,\quad\text{ then one obtains }\quad 0<\inf \rho_{g,h}\leq \sup \rho_{g,h}<\infty,
$$
so there is the trivial bounded identification operator
$$
I_{g,h}:L^2(M,g)\longrightarrow L^2(M,h),\quad I_{g,h}f:= f.
$$
Note that, however, the adjoint of this operator is nontrivial multiplication operators.\vspace{1mm}

The following result is the main result of this note, and it generalizes Corollary A from \cite{guneysu_scattering_2020} from manifolds with Ricci curvature bounded from below by a constant to Kato-Ricci manifolds:

\begin{theorem}\label{main}
Assume that $g,h$ are complete and quasi-isometric Riemannian metrics on $M$ such that $(M,g)$ and $(M,h)$ are Kato-Ricci with
\begin{align}\label{assu}
\int_X \mu_l(x,1)^{-1} \delta_{g,h}(x)\,d\mu_l(x) < \infty\quad\text{for some/both $l=g,h$.}
\end{align}
Then the wave operators $W_\pm(-\Delta_h, -\Delta_g;I_{g,h})$ exist and are complete. 
\end{theorem}

\begin{proof} Clearly, by quasi-isometry one has (\ref{assu}) for some metric, if and only if this holds for both metrics. By the proof of Theorem A in \cite{guneysu_scattering_2020}, it suffices to show that the operators 
\begin{align}\label{drei}
\delta_{g,h}P^g_t, \quad \delta_{g,h}\hat{P}^g_t ,\quad \delta_{g,h}P^h_t, \quad \delta_{g,h}\hat{P}^h_t 
\end{align}
are Hilbert-Schmidt for some $t>0$. To this end, it suffices to show that for every Kato-Ricci manifold $M$ (omitting the depencence of the metric in the notation) one has the following estimates: for all $0<t<1$, $f\in L^2(M)$, $x\in M$,
\begin{align}
\label{null}|P_tf(x)| \lesssim  \mu(x,t)^{-1/2}\left(\int_M |f(y)|^2 d\mu(y)\right)^{1/2},\\
\label{eins}|\hat{P}_t f(x)| \lesssim t^{-1/2}\mu(x,t)^{-1/2} \left(\int_M |f(y)|^2 d\mu(y)\right)^{1/2}.
\end{align}
Indeed, then Riesz-Fischer's duality theorem implies that the integral kernel $[\delta_{g,h}P^g_t](x,y)$ of $\delta_{g,h}P^g_t$ satisfies
$$
\int_M |[\delta_{g,h}P^g_t](x,y)|^2d\mu_g(y)\lesssim  \delta_{g,h}(x)\mu_g(x,t)^{-1}\quad\text{for all $0<t<1$, $x\in M$,} 
$$
and so
$$
\int_M |[\delta_{g,h}P^g_t](x,y)|^2(x,y)d\mu_g(y)d\mu_g(y)\mid_{t=1}<\infty
$$
by (\ref{assu}), and similiarly for the other three operators in (\ref{drei}).\vspace{1mm}

To prove (\ref{null}), we record that (cf. Theorem 3.3 in \cite{guneysu_locally_2023}) as a consequence of the Kato-Ricci assumption, there exist constants $\alpha, \beta, \gamma>0$ such that for all $t>0$ and $x,y\in M$,
\begin{align}\label{GUE}
P(t, x, y) \leq \alpha \mu(x, \sqrt{t})^{-1} e^{-\frac{\beta d(x,y)^2}{t}} e^{\gamma t}.
\end{align}
Thus for $0<t<1$, using Cauchy-Schwarz and (\ref{assq}),
\begin{align*}
|P_tf(x)|&\leq \int_M P(t, x, y) |f(y)|d\mu(y)\\
&\leq \left(\int_M P(t, x, y)^2d\mu(y)\right)^{1/2} \left(\int_M|f(y)|^2d\mu(y)\right)^{1/2} \\
& \lesssim \mu(x,t)^{-1/2}\left(\int_M P(t, x, y)d\mu(y)\right)^{1/2} \left(\int_M|f(y)|^2d\mu(y)\right)^{1/2}\\
&\leq \mu(x,t)^{-1/2} \left(\int_M|f(y)|^2d\mu(y)\right)^{1/2}.
\end{align*}
To prove (\ref{eins}), we use that by Theorem 1.5 in \cite{braun_heat_2021}, for every $t>0$, $f\in  L^2(M)$, $x\in M$, $\xi\in T_xM$, we have the Bismut-Elworthy-Li formula
\[
( \hat{P}_tf(x), \xi ) = \frac{1}{t} \E \left[ f(X_s(x)) \int_0^t \big( Q_s(x)\xi, dW_s(x) \big) \right],
\]
where
\begin{itemize}
\item $X(x)$ is a Brownian motion\footnote{We understand Brownian motion as a $\Delta$ diffusion, rather than a $\Delta/2$ diffusion.} in $M$ starting from $x$
\item with $\partr^x_s : T_xM \to T_{X_s(x)}M$, $s\geq 0$, the stochastic parallel transport w.r.t. the Levi-Civita connection along the paths of $X(x)$, the $\End(T_xM)$-valued process $Q(x)$ is defined as the unique solution to the pathwise ordinary differential equation
\begin{align}\label{gron}
dQ^x_s = -  Q^x_s (\partr^x_s)^{-1}\ric(X_s(x)) \partr_s^x ds,\quad Q^x_0 = \id_{T_xM};
\end{align}
here $\ric(X_s(x))$ is interpreted as an endomorphism of $T_{X_s(x)}M$
\item $W(x)$ is a Brownian motion on the Euclidean space $T_xM$
\item $\int_0^\bullet \big( Q_s(x)\xi, dW_s(x) \big)$ denotes the Itô integral.
\end{itemize}
By Cauchy-Schwarz, 
\[
|(  \hat{P}_tf(x),\xi )| \leq \frac{1}{t} \E \left[ |f(X_s(x))|^2 \right]^{1/2}
\E \left[\left|\int_0^t \langle Q_s^x\xi, dW_r^x \rangle\right|^2\right]^{1/2},
\]
so the first term can be estimated using that the transition density of Brownian motion is given by the heat kernel and (\ref{GUE}), 
\[
\E \left[|f(X_s(x))^2\right] = \int_M |f(y)|^2 p(t,x,y)\,d\mu(y)  \lesssim \mu(x,\sqrt{t})^{-1} \int_M |f(y)|^2\,d\mu(y) .
\]
For the second term, we use the Itô isometry and 
$$
|Q^x_s\zeta |\leq e^{\int_0^s \sigma_g^-(X_s(x))\,ds} \quad\text{ for all $|\xi|\leq 1$, $s\geq 0$,}
$$
a consequence of (\ref{gron}) and Gronwall's lemma, to deduce for all $0<t<1$,
\[
\E\left[ \left| \int_0^t \langle Q^x_s\xi, dW_s^x \rangle \right|^2 \right]
= 2\E\left[\int_0^t |Q_s^x\zeta |^2\,ds\right]  \lesssim t \E\left[ e^{2\int_0^1 \sigma_g^-(X_s(x))\,ds} \right],
\]
which is finite uniformly in $x$ by the Kato-Ricci assumption and Khas'minskiii's lemma (Lemma VI.8 in \cite{bookbatu}). This completes the proof of (\ref{eins}).

\end{proof}

In the context of geometric scattering theory, there is also another natural identification operator other than the above $I_{g,h}$, which is well-defined even if the metrics are not quasi-isometric, namely, the unitary operator
$$
J_{g,h}:L^2(X,\mu_g)\longrightarrow L^2(X,\mu_h),\quad J_{g,h}\psi:= \psi/\sqrt{\rho_{g,h}},
$$
which induces the wave operators $W_\pm(-\Delta_h, -\Delta_g;J_{g,h})$. A criterion which makes sure that these two choices of wave operators are equal has been given in \cite{asymptotic}. In fact, the following result follows immediately from (\ref{null}) and the proof of Theorem 4.1 in \cite{asymptotic}:

\begin{theorem} Under the assumptions of Theorem \ref{main}, one has
$$
W_\pm(-\Delta_h, -\Delta_g;J_{g,h})=W_\pm(-\Delta_h, -\Delta_g;I_{g,h}).
$$
\end{theorem}

\bibliography{biblio.bib}
\bibliographystyle{plain}
\end{document}